\documentclass[review]{elsarticle}
\usepackage{hyperref}
\usepackage{amsfonts}

\begin{document}

\newtheorem{thm}{Theorem}
\newtheorem{lem}{Lemma}
\newtheorem{cor}{Corollary}
\newtheorem{pro}{Proposition}
\newtheorem{conj}{Conjecture}
\newdefinition{dfn}{Definition}
\newproof{prf}{Proof}


\title{A Semigroup Proof of the Bounded Degree Case of S.B. Rao's Conjecture
on Degree Sequences and a Bipartite Analogue}
\author[cja]{Christian Joseph Altomare}
\ead{altomare@math.ohio-state.edu}
\address{The Ohio State University, 231, West 18th Avenue,
Columbus, Ohio, United States}

\begin{abstract}
S.B. Rao conjectured in 1971 that graphic degree sequences are well
quasi ordered by a relation $\preceq$ defined in terms of the induced
subgraph relation\cite{rao_private}.
In 2008, M. Chudnovsky and P. Seymour proved this
long standing Rao's Conjecture by giving structure theorems for graphic
degree sequences\cite{proof_rao}.

In this paper, we prove and use a variant of Dickson's
Lemma\cite{dicksons_lemma} from
commutative semigroup theory to give
a short proof of the bounded degree case of Rao's Conjecture that is
independent of the Chudnovsky-Seymour structure theory. In fact, we
affirmatively answer two questions of N. Robertson\cite{neil_private},
the first of which implies the bounded degree case of Rao's Conjecture.
\end{abstract}

\maketitle




\section{Introduction}


Let $G$ be a finite, simple graph and let $D(G)=(d_1,\ldots,d_n)$
be its list of vertex degrees listed in decreasing order.
The sequence $D(G)$
is known as the degree sequence of $G$, and $G$ is said to realize
$D$. A sequence $(d_1,\ldots,d_n)$ of nonnegative integers is said
to be a graphic degree sequence if it is realized by some graph.
Given graphic degree sequences $D_1$ and $D_2$, we define $D_1\preceq D_2$
to mean there is $G_1$ realizing $D_1$ and $G_2$ realizing $D_2$
such that $G_1\sqsubseteq G_2$, where $\sqsubseteq$ is the induced
subgraph relation. The reader may check
that $\preceq$ is a transitive relation on degree sequences.
For other basic graph theoretic definitions, we refer the
reader to \cite{diestel_book}.

We recall that a quasi order $(Q,\le)$ is a reflexive, transitive
relation $\le$ on a class~$Q$.
A quasi order $(Q,\le)$ is said to be a well quasi order
if $Q$ contains no infinite decreasing sequence and no infinite
antichain. Equivalently, $(Q,\le)$ is a well quasi order if for every
infinite sequence $q_1,q_2,\ldots$
in $Q$ there are positive
integers $i<j$ such that $q_i\le q_j$.

With these definitions, we may state Rao's Conjecture, posed in
1971 by S.B. Rao\cite{rao_private}
and finally proved in 2008 by M. Chudnovsky and
P. Seymour\cite{proof_rao}.

\def\mod{\hbox{mod}}

\begin{thm} \label{raos_conj}
Degree sequences of finite graphs are well quasi ordered by $\preceq$.
\end{thm}

Independently, N. Robertson had asked\cite{neil_private} if
graphic degree sequences of bounded degree can be realized
as disjoint unions of graphs with bounded sized components,
noticing that an affirmative answer would imply the
bounded degree case of \ref{raos_conj}.
Motivated by this question, he further asked
for a bipartite analogue. Namely, Robertson asked if degree
sequences of bipartite graphs of bounded degree can be realized
as disjoint unions of bipartite graphs with bounded sized
components\cite{neil_private}.

In this work, we use a variant of Dickson's Lemma from
commutative semigroup theory to prove a more general fact
that yields affirmative answers to both of Robertson's questions
as corollaries. In particular, we obtain a new proof of
the bounded degree case of \ref{raos_conj} that does
not depend on the Chudnovsky-Seymour structure theory.
This proof abstracts a proof given in the author's doctoral
thesis\cite{my_thesis}.

While our proof has the disadvantage of only going through for bounded
degree, it is fairly short. Moreover,
our proof is no longer restricted to
graphs and goes through equally well for partial orders,
hypergraphs, or any class of structured sets at all for which
nonnegative integers can be assigned to each point in a way
that respect disjoint
union and such that regular elements exist.
In particular, even for graphs, these nonnegative integers
need no longer represent the degree.
This is worth noting since some of the most commonly
used tools for degree sequences, such as switchings\cite{switching} and the
Erd\"os-Gallai inequalities\cite{erdos_gallai},
have no known counterparts in this
more general setting.

\section{The Semigroup Lemma}

A commutative semigroup is a set $S$ together with an
associative, commutative binary operation $+$. We need not assume
existence of an identity element. For basic facts and terminology,
we refer the reader to \cite{grillet}, but our presentation is
self contained.
Given a semigroup $(S,+)$ and
subsets $Y\subseteq X$ of $S$, we say that $Y$ generates $X$ if
every $x$ in $X$ can be written as $y_1+y_2+\cdots+y_n$ for
some points $y_1,\ldots,y_n$ of $Y$. We say that $X$ is finitely
generated if some finite subset $Y$  of $X$ generates $X$.

We now work exclusively with the free commutative semigroup
$\mathbb{N}^k$, where we assume $k$ is fixed throughout.
Given $x=(x_1,\ldots,x_k)$ in $\mathbb{N}^k$, the support $\hbox{supp}(x)$
is defined as the set of $i$ such that $x_i>0$.

\begin{dfn}
Let $X$ be a subset of $\mathbb{N}^k$. We say that $X$ is grounded if for all
$i$ in $\{1,\ldots,k\}$ there is $x$ in $X$ with $\hbox{supp}(x)=\{i\}$.
\end{dfn}

Let $(x_1,\ldots,x_k)$, $(y_1,\ldots,y_k)$,
and $(t_1,\ldots,t_k)$ be elements in $\mathbb{N}^k$. We say that
$$(x_1,\ldots,x_k)\equiv(y_1,\ldots,y_k)\,\mod\, (t_1,\ldots,t_k)$$
if $x_i\equiv y_i\,\mod\, t_i$ for each $i$.

\begin{lem}\label{semlem} 
Every grounded subset of $(\mathbb{N}^k,+)$ is finitely generated.
\end{lem}

\begin{prf}
Fix a grounded set $X$. Then for each $i$ in $\{1,\ldots,k\}$,
we may choose an element of the form $(0,\ldots,0,t_i,0,\ldots,0)$
in $X$,
where $t_i>0$ occurs in position~$i$. Without loss of generality,
we may choose the minimal such $t_i$ for each $i$.
Note that equivalence
modulo $(t_1,\ldots,t_k)$ is an equivalence relation $\sim$ on $\mathbb{N}^k$
with only finitely many equivalence classes.

The partial order $(\mathbb{N},\le)$ with the usual ordering of the
natural numbers is obviously a well quasi order.
Since the
product of finitely many well quasi orders is a well quasi order,
we see that $(\mathbb{N}^k,\le)$ is well quasi ordered  when considered as
a product order. In particular, every antichain in $(\mathbb{N}^k,\le)$
is finite.

Given a nonempty $\sim$ class $C$, the (possibly empty) set $M_C$ of
$(\mathbb{N}^k,\le)-\{(0,\ldots,0)\}$
minimal elements of $C$ is an antichain in $(\mathbb{N}^k,\le)$
and therefore finite by
the previous paragraph. Let $Y$ be the union over all
$\sim$ classes $C$ of the sets $M_C$.
Then $Y$ is the finite union of finite sets
and so is finite. It is thus enough to show
$Y'=Y\cup\{(0,\ldots,0)\}$ generates $X$.

Choose $x=(x_1,\ldots,x_k)$ in $X$.
If $x=(0,\ldots,0)$ then $x$ is in $Y'$ and we are done. Assume not.
Let $C$ be the $\sim$ class of $x$.
Since $(\mathbb{N}^k,\le)$ is well founded, $C$ contains a
$(\mathbb{N}^k,\le)-\{(0,\ldots,0)\}$ minimal element
$(m_1,\ldots,m_k)$ such that $(m_1,\ldots,m_k)\le(x_1,\ldots,x_k)$.
Then $m_i\le x_i$ for each $i$. Since
$(m_1,\ldots,m_k)\sim (x_1,\ldots,x_k)$ by hypothesis, we see
that for each $i$, the equation $x_i-m_i=c_it_i$ holds
for some nonnegative integer $c_i$.
Therefore
$$(x_1,\ldots,x_k)=
(m_1,\ldots,m_k)+\sum_{i=1}^k c_i(0,\ldots,0,t_i,0,\ldots,0).$$
We know $(m_1,\ldots,m_k)$ is in $Y'$ by hypothesis. It is easy
to see that $(0,\ldots,0,t_i,0,\ldots,0)$ is a minimal nonzero
element in its $\sim$
class.
Therefore $(0,\ldots,0,t_i,0,\ldots,0)$ is
in $Y'$ as well, by which we see $Y'$ generates
$(x_1,\ldots,x_k)$. As $(x_1,\ldots,x_k)$ in $X$ was chosen arbitrarily,
we see that $Y'$ generates $X$ as claimed.
\end{prf}

\section{Structured Sets}

Our main theorem will apply equally to the class of finite graphs
and the class of finite, bipartite graphs, the class of finite
posets, and so on, so we need a general way to speak of all these
classes of objects. As numerous mathematical objects are defined
as a set together with some structure on it, which could be a binary
operation, a relation, a set of subsets, and so on, rather than
try to define some generalized structure on a set that includes all
these things, we prefer to consider the structure as nothing more than
a label. For instance, a partial order $(P,\le)$ would be considered
a set $P$ together with label $\le$. For us then, a structured set
is simply a set together with a label.

The one operation we need on our structured sets is that of coproduct,
which will correspond to $+$ in the semigroup.
For our purposes,
we simply consider $\coprod$ as an arbitrary associative, commutative
binary operation on a
class of structured sets. We make this and the previous paragraph
precise in the following definition.

\begin{dfn}
A structured set class is a class $U$ of ordered pairs
together with an associative, commutative binary operation
$\coprod:U\times U\to U$
called coproduct such that
$P$ is a finite set for each ordered pair $(P,T)$ in $U$.
We call members of $U$ structured sets. 
\end{dfn}

When no confusion arises,
we sometimes say $P$ instead of $(P,T)$.

Note that in the above definition, $\coprod$ is not a
function in the sense of being a set of ordered pairs. The
binary operation $\coprod$ is, in natural cases, a proper class of
ordered pairs as $U$ is. While this fact is worth noting,
it creates no problems, and we do not concern ourselves
with such foundational issues here. We use proper classes
freely and without comment.

Also note that since $\coprod$ is both associative and commutative,
we could consider $\coprod$ itself as a semigroup whose domain is
a proper class. Though formally correct, we do not take this point
of view, as we find it is a greater
aid to intuition to think of $\coprod$ as a coproduct of structures
than as addition in an abelian semigroup.

We have generalized the notion of a class of finite graphs
to the notion of a structured set class. We now need to
generalize the notion of the degree sequence of a graph to
this new setting. In fact, using degree sequences for graphs
would not allow us to use the semigroup lemma as even graphs
of bounded degree may have arbitrarily long degree sequences.
The degree sequences of graphs of degree at most $k$ are not,
therefore, contained in $\mathbb{N}^r$ for any $r$.

The solution is to instead use what we define as the regularity
sequence of a graph. Given a finite graph $G$, its regularity
sequence is the unique sequence $R_G$ such that
$$R_G(i)=|\{v:v\in G\hbox{ and }d_G(v)=i\}|$$
for each natural number $i$. It is simple to check that the degree
sequence and regularity sequence of a graph each uniquely
determine the other. Instead of generalizing the notion of
degree sequence to structured set classes, we generalize the
notion of regularity sequences to structured set classes.

\begin{dfn}
Let $U$ be a structured set class. 
A structured set function for the class $U$ is a function $F$
whose range is a subset of $\mathbb{N}$ and
whose domain consists of all triples $(P,T,x)$ such that
$(P,T)$ is a structured set in $U$ and $x$ is in $P$.
\end{dfn}

\begin{dfn}
Let $U$ be a structured set class.
Let $F$ be a structured set function for $U$ and
let $(P,T)$ in $U$ be a structured set.
The $F$-regularity sequence $R_{F,P}$ of $(P,T)$
is defined by letting
$$R_{F,P}(i)=|\{v:v\in P\hbox{ and }F(P,T,v)=i\}|$$
 for each natural number $i$.
\end{dfn}

Note that a regularity sequence is, in particular, a sequence.
We may therefore add regularity sequences.

\begin{dfn}
Let $U$ be a structured set class.
A structured set function $F$ for $U$ is called additive if for all
structured sets $P$ and $Q$, we have
$$R_{F,P\coprod Q}=R_{F,P}+R_{F,Q}.$$
\end{dfn}

The reader may check that if $U$ is the class of finite graphs,
considered as structured sets by letting
the edge set $E(G)$ be the label of
the finite vertex set $V(G)$ and $\coprod$ the disjoint union of graphs,
then the structured set function $F$ taking a triple
$(V(G),E(G),v)$ to $d_G(v)$ is additive.
Additivity
of this $F$ simply states that the number of vertices of degree $i$
in the disjoint union of $G$ and $H$ is the number of vertices of
degree $i$ in $G$ plus the number of vertices of degree $i$ in $H$.

The following definition generalizes to our
new setting the notion of two graphs having the same 
regularity sequence.

\begin{dfn}
Let $U$ be a structured set class.
Let $F$ be a structured set function for $U$. Let $P$ and $Q$ be
structured sets. We say that
$P$ and $Q$ are $F$-equivalent if $R_{F,P}=R_{F,Q}$.
\end{dfn}

\begin{dfn}
Let $U$ be a structured set class and $F$ be a
structured set function for $U$. We say that $U$
is $F$-finitely representable if there is some finite subset $Z$ of $U$
such that every structured set in $U$ is $F$-equivalent to some
structured set of the form
$$\coprod_{i=1}^n P_i,$$
with each $P_i$ in $Z$.
\end{dfn}

The following definition abstracts the notion of the class of finite
graphs with all degrees at most~$k$.

\begin{dfn}
Let $U$ be a structured set class, $F$ a structured set function for $U$,
and $k$ a nonnegative integer. Then $U_{F,k}$ denotes the class of
structured sets
$(P,T)$ such that $R_{F,P}(i)=0$ for all integers $i>k$.
\end{dfn}

We now make a definition that
will be used as a hypothesis to ensure a subset of a semigroup is
grounded.

\begin{dfn}
Let $U$ be a structured set class and $F$ a structured set function for $U$.
Then $F$ is said to have regulars if for all nonnegative integers $i$ there
is a structured set $P$ such that $R_{F,P}$ has support $\{i\}$.
\end{dfn}

Note that for the structured set class~$U$ of finite graphs and the structured
set function~$F$ taking each vertex to its degree, $F$ has regulars since
there are $k$-regular graphs for each nonnegative integer~$k$.

\section{The Main Theorems}

We now state our main theorem.

\begin{thm}\label{main} 
Let $U$ be a structured set class,
$k$ a nonnegative integer, and $F$ an additive
structured set function for $U$ that has regulars.
Then $U_{F,k}$ is $F$-finitely representable.
\end{thm}


\begin{prf}
By definition
of $U_{F,k}$, we know that for each structured set $P$ in $U_{F,k}$ and $i>k$
that $R_{F,P}(i)=0$. We may therefore think of the regularity
sequence $R_{F,P}$ as the finite sequence
$R_{F,P}(0),\ldots,R_{F,P}(k)$ of length $k+1$, which we consider
as an element of the additive semigroup $(\mathbb{N}^{k+1},+)$. Let
$X$ be the set of points in $(\mathbb{N}^{k+1},+)$ corresponding to
regularity sequences of
structured sets in $U_{F,k}$.

Since $F$ has regulars, we see that $X$ is grounded.
By \ref{semlem}, we know that $X$ is finitely generated. Let $Y$
be a finite generating set. Each member $y$ of $Y$ is in particular a
member of $X$, and therefore there is a structured set $P$ in $U_{F,k}$ that
has $y$ as its regularity sequence. We may thus choose a finite set $Z$ of
structured sets in~$U_{F,k}$ such that each regularity sequence $y$ in $Y$ is the
regularity sequence of some structured set in~$Z$.

Now take an arbitrary structured set $P$ in $U_{F,k}$.
We know its regularity sequence
$x$ is in $X$ by definition. Therefore
$$x=a_1y_1+\ldots+a_ny_n$$
for some $n\ge 1$, nonnegative integers $a_i$, and members $y_i$ of $Y$.
Therefore $R_{F,P}=R_{F,Q}$, where $Q$ is the structured set
$$\coprod_{i=1}^n a_iQ_i,$$
where $Q_i$ is a structured set in $Z$ with regularity sequence $y_i$
and $a_iQ_i$ denotes the structured set
$$\coprod_{j=1}^{a_i}Q_i.$$

We therefore see there is a finite set $Z$ of structured sets in $U_{F,k}$
such that each structured set in $U_{F,k}$ is $F$-equivalent to a coproduct
of structured sets in $Z$. By definition of $F$-finite representability, this
completes the proof.
\end{prf}

We now apply this theorem to answer Robertson's original questions. Though
Robertson asked if graphic degree sequences of bounded degree may be
realized with bounded sized components, we note this is equivalent to asking
if graphic degree sequences of bounded degree may be realized as disjoint
unions of graphs from a fixed finite set, and similarly for the bipartite
analogue. We find this reformulation somewhat more convenient as then
\ref{main} more directly applies.

\begin{cor}\label{cor1}
Degree sequences of finite, bipartite graphs with bounded degree can be
realized as disjoint unions of bipartite graphs from a fixed finite set.
\end{cor}

\begin{prf}
Let $U$ be the class of finite, bipartite graphs $G=(V(G),E(G))$,
considered as a structured set class by choosing label $E(G)$ for the
set $V(G)$,
$\coprod$ representing disjoint union,
and $F$ the additive structured set function for $U$
taking $(V(G),E(G),v)$ to $d_G(v)$.
To show that $F$ has regulars,
simply note that $K_{j,j}$ is a $j$-regular bipartite graph for
each nonnegative $j$. We thus see by \ref{main} that
regularity sequences of finite, bipartite graphs with bounded degree can be
realized as disjoint unions of bipartite graphs from a fixed finite set~$Z$.
It is trivial that the same is thus true for degree sequences.
\end{prf}

\begin{cor}\label{cor2}
Degree sequences of finite graphs with bounded degree can be
realized as disjoint unions of graphs from a fixed finite set.
\end{cor}

\begin{prf}
The proof is exactly as in that of the previous corollary except we let
$U$ be the structured set class of finite graphs.
\end{prf}

It is worth noting that neither \ref{cor1} nor \ref{cor2} is stronger
than the other, as both the hypotheses and the conclusions of \ref{cor1}
are stronger than that of \ref{cor2}. We now give the simple proof that
\ref{cor2} implies the bounded degree case of Rao's Conjecture.

\begin{cor}
Fix $k$. Degree sequences of finite graphs with degrees at most~$k$
are well quasi ordered by $\preceq$.
\end{cor}

\begin{prf}
By \ref{cor2}, there is a finite set $Z$ of finite graphs with degrees
at most $k$ such that every graphic degree sequence with degrees at
most~$k$ can be realized as a disjoint union of graphs in $Z$.
Since degree sequences and regularity sequences contain the same information,
we may consider $\preceq$ as a relation on regularity sequences, and we
again think of regularity sequences as points in $(\mathbb{N}^{k+1},+)$.

We note that given points $x$ and $x'$ in $\mathbb{N}^{k+1}$, if
$x\le x'$ in the product order $(\mathbb{N}^{k+1},\le)$ then
$x\preceq x'$. This implies every $\preceq$ antichain is a $\le$
antichain. We have previously noted that all $\le$ antichains are
finite, which implies all $\preceq$ antichains are finite, thus
completing the proof.
\end{prf}


\bibliographystyle{plain}
\bibliography{mybib}
\end{document}